\documentclass[10pt]{amsart}  
\usepackage{latexsym}
\usepackage{amsmath}
\usepackage{amssymb}
\usepackage{mathrsfs}
\usepackage{graphicx}
\usepackage{bbold}
\input epsf
\input xy
\xyoption{all}
\input diagxy

\def\st{\mathop{\fam0 st}\nolimits}
\def\St{\mathop{\fam0 St}\nolimits}

\def\fin{\mathop{\fam0 fin}\nolimits}
\def\fil{\mathop{\fam0 fil}\nolimits}
\def\mix{\mathop{\fam0 fil}\nolimits}

\def\Orth{\mathop{\fam0 Orth}\nolimits}

\def\dom{\mathop{\fam0 dom}\nolimits}
\def\On{\mathop{\fam0 On}\nolimits}

\def\ZFC{\mathop{\fam0 ZFC}}
\def\widevert{\kern1pt \vrule height7.5pt depth2.5pt width1pt\kern1pt}
\def\Rightharpoonup{\mathrel{\vcenter{\offinterlineskip\halign{\hfil##\hfil\cr
$\rightharpoonup$\cr\hskip.6pt\hrulefill\hskip.6pt\cr\noalign{\vskip.5pt}}}}}

\def\upwardarrow{\mathord{
  \hbox to 5pt{\hss$\vcenter{\hbox to 2.4pt{\hss$\mathchar"222$\hss}\hrule}\hss$}
}}
\def\downwardarrow{\mathord{
  \hbox to 5pt{\hss$\vcenter{\hrule\hbox to 2.4pt{\hss$\mathchar"223$\hss}}\hss$}
}}

\begin{document}
\def\?{?\vadjust{\vbox to 0pt{\vss\hbox{\kern\hsize\kern1em\bf ?!}}}}

\title
[Leibnizian, Robinsonian, and Boolean Valued Monads]
{Leibnizian, Robinsonian,\\ and Boolean Valued Monads}

\author{S.~S. Kutateladze}
\address[]{
Sobolev Institute of Mathematics\newline
\indent 4 Koptyug Avenue\newline
\indent Novosibirsk, 630090\newline
\indent Russia}
\email{
sskut@math.nsc.ru
}
\begin{abstract}
This is an overview of  the present-day
versions of monadology with some applications
to vector lattices and linear inequalities.
Two approaches to combining nonstandard set-theoretic models are 
sketched and illustrated  by order convergence, 
principal projection, and polyhedrality.

\end{abstract}

%\date{June 6, 2011}
%\thanks{This article bases on a talk at
%the 20th St.~Petersburg Summer Meeting in Mathematical Analysis,
%June 24--29, 2011.}

\maketitle

The notion of monad is central to  external set theory.
Justifying  the simultaneous use of infinitesimals
and the technique of descending and
ascending in vector lattice theory
requires adaptation of monadology for the
implementation of filters in  Boolean valued universes.
This is still a rather uncharted area of research.
The two approaches are available now.
One is to apply monadology  to the descents of objects.
The other consists in
applying the standard monadology inside the Boolean valued universe ${\mathbb V}^{(\mathbb B)}$
over a complete Boolean algebra $\mathbb B$, while ascending and descending by
the Escher rules (cp.~\cite{IBA} and \cite{BooPos}).

These approaches are sketched and illustrated by tests for order convergence
and rules for fragmenting and projecting positive operators in vector lattices.
Also,
Lagrange's principle is shortly addressed  in  polyhedral environment with inexact data.

\section*{Basics of Monadology}
The concept of monad stems from Ancient Greece.
Monadology as a~philosophical  doctrine
is a creation of Leibniz (cp.~\cite{Leibniz} and \cite{MVL}).
The general theory of the monads of filters was proposed by
Luxemburg (cp.~\cite{Luxemburg}) within  Robinson's
nonstandard analysis (cp.~\cite{Dauben}).

Let ${\mathscr F}$ be a~standard filter; ${} ^\circ{\mathscr F}$,
the standard core of ${\mathscr F}$; and
${}^a\!{\mathscr F}\!:= {\mathscr F}\setminus{}^\circ\!{\mathscr F}$,
the external set of {\it remote\/}
elements
of $\mathscr F$. Note that
$$
\mu ({\mathscr F})\!:= \bigcap {}^\circ\!{\mathscr F} =\bigcup {}^a\!{\mathscr F}
$$
is the  {\it  monad\/}
of ${\mathscr F}$.
Also,
${\mathscr F}={}^\ast \fil{(\{ \mu({\mathscr F}) \})}$; i.e., $\mathscr F$
is the standardization of the collection
$\fil{(\mu({\mathscr F}))}$
of all supersets of $\mu({\mathscr F})$.

Let
$\mathscr  A$
be a~filter on
$X\times Y$,
and let
$\mathscr  B$
be a~filter on
$Y\times Z$.
Put
$
\mathscr  B\circ \mathscr  A
:=\fil
\{B\circ A\mid A\in\mathscr  A,\, B\in \mathscr  B\},
$
where we may assume all $B\circ A$
nonempty. Then
$$
\mu(\mathscr  B\circ\mathscr  A)=
\mu(\mathscr  B)\circ\mu(\mathscr  A).
$$

{\bf Granted Horizon Principle.}
{\it Let
$X$
and
$Y$
be standard sets. Assume further that
$\mathscr F $
and
$\mathscr  G $
are standard filters on~$X$ and~$Y$ respectively
satisfying
$\mu (\mathscr  F )\cap {}^\circ \!X\ne\varnothing $.
Distinguish  a~remote set $F$ in ${}^a\!\mathscr  F$.
Given a~standard correspondence
$f\subset X\times Y$
meeting $\mathscr  F $,
the following  are equivalent:

{\bf (1)}
$f(\mu (\mathscr  F )-F)\subset\mu (\mathscr  G )$;

{\bf (2)}
$(\forall \,F'\in {}^a\!\mathscr  F )\, f(F'-F)\subset\mu (\mathscr  G )$;

{\bf (3)}
$f(\mu (\mathscr  F ))\subset\mu (\mathscr  G )$.
}

\section*{Filters within ${\mathbb V}^{(\mathbb B)}$}

Let
${\mathbb B}$
be a~complete Boolean algebra. Given an ordinal
$\alpha$,
put
$$
V_{\alpha}^{({\mathbb B})}
:=\{x \mid
(\exists \beta\in\alpha)\ x:\dom (x)\rightarrow
{\mathbb B}, \dom (x)\subset V_{\beta}^{({\mathbb B})}  \}.
$$

The {\it Boolean valued
universe\/}
${\mathbb V}^{({\mathbb B})}$
is
$$
{\mathbb V}^{({\mathbb B})}:=\bigcup\limits_{\alpha\in\On} V_{\alpha}^{({\mathbb B})},
$$
with $\On$ the class of all ordinals.

The truth
value $[\![\varphi]\!]\in {\mathbb B}$ is assigned to each formula
$\varphi$ of $\ZFC$ relativized to ${\mathbb V}^{({\mathbb B})}$.

Let $Q$ be the Stone space of a~complete Boolean algebra~$\mathbb B$.
Denote by $\mathfrak U$  the (separated) Boolean valued
universe~${\mathbb V}^{(\mathbb B)}$.
Given $q\in Q$, put
 $
 u\sim_q v\ \leftrightarrow\ q\in[\![u=v]\!].
 $
Consider the bundle
 $$
 V^Q:=\big\{\big(q,{\sim_q}(u)\big) \mid q\in Q,\, u\in\mathfrak{U}\big\}
 $$
 and denote  $\big(q,{\sim}_q(u)\big)$ by $\widehat u(q)$.
Hence $\widehat u:q\mapsto \widehat u(q)$ is a section of~$V^Q$ for every $u\in\mathfrak{U}$.
 Note that to each $x\in V^Q$ there are $u\in\mathfrak{U}$
 and $q\in Q$ satisfying $\widehat u(q)=x$.
Moreover, we have $\widehat u(q)=\widehat v(q)$
 if and only if $q\in[\![u=v]\!]$.

Make each fiber~$V^q$ of~$V^Q$ into an algebraic system of
 signature~$\{\in\}$ by letting
 $
 V^q\models x\,{\in}\, y \ \leftrightarrow\  q\in[\![u\,{\in}\, v]\!],
 $
where $u,v\in\mathfrak{U}$ are such that $\widehat u(q)=x$ and $\widehat v(q)=y$.

The class  $\{\widehat u(A)\mid u\in\mathfrak U\}$,
 with $A$ a~clopen subset of~$Q$, is a base for some topology
 on~$V^Q$. Thus~$V^Q$
 as a~continuous bundle called  a~{\it continuous polyverse}.
By a~{\it continuous section\/} of~$V^Q$  we mean a section
that is a continuous function.  Denote by $\mathfrak C$
 the class of all continuous sections of~$V^Q$.

 The mapping $u\mapsto\widehat u$ is a bijection between~$\mathfrak U$
 and $\mathfrak C$, yielding a convenient functional realization of the
 Boolean valued universe ${\mathbb V}^{(\mathbb B)}$.
 This universal construction belongs to
 Gutman and Losenkov (cp.~\cite{GutLos}).

 The functional realization visualizes descending and ascending, the Escher rules,
 and the  Gordon Theorem (cp.~\cite{Trends}).

Let ${\mathscr G}$ be a~filterbase on $X$, with
$X\in{\mathscr P}({\mathbb V}^{(\mathbb B)})$.
Put
$$
\gathered
{\mathscr G}^\prime\!:= \{ F \in{\mathscr P(X\!\!\uparrow)\!\!\downarrow}\mid
(\exists G \in {\mathscr G})\,[\![\,F \supset {G\!\!\uparrow}\,]\!]={\mathbb 1}\};\\
{\mathscr G}^{\prime\prime}\!:=\{{G\!\!\uparrow}\mid G\in{\mathscr G}\}.
\endgathered
$$
 Then
${{\mathscr G}^\prime\!\!\uparrow}$
and
${{\mathscr G}^{\prime \prime}\!\!\uparrow}$
are bases of  the same filter
${\mathscr G}^\uparrow$
on
${X\!\!\uparrow}$
inside
${\mathbb V}^{(\mathbb B)}$---the {\it ascent\/}
of
${\mathscr G}$.
If $\mix({\mathscr G})$
is the set of all mixings of nonempty families of elements of
${\mathscr G}$ and ${\mathscr G}$ consists of cyclic sets; then
$\mix({\mathscr G}) $ is a~filterbase on $X$ and
${\mathscr G}^\uparrow=\mix({\mathscr G})^\uparrow$.

If ${\mathscr F}$ is a~filter on $X$ inside
${\mathbb V}^{(\mathbb B)}$
then  put
${\mathscr F}^\downarrow\!:= \fil{({\{{F\!\!\downarrow} \mid
F\in {\mathscr F}\!\!\downarrow\}})}$.
The filter
${\mathscr F}^\downarrow$
is the {\it  descent\/}
of ${\mathscr F}$.
A filterbase
${\mathscr G}$ on ${X\!\!\downarrow}$ is {\it extensional\/}
provided that
$\fil{({\mathscr G})}={\mathscr F}$
for some filter ${\mathscr F}$ on $X$.

The descent of an ultrafilter on $X$ is
a {\it proultrafilter\/}
on ${X\!\!\downarrow}$. A~filter with a~base of cyclic sets is
{\it cyclic}.
Proultrafilters are maximal cyclic filters.

Fix a~standard complete Boolean algebra $\mathbb B$ and
think of $ {\mathbb V}^{(\mathbb B)}$
to be  composed of internal sets. If $A$ is  external
then the {\it cyclic hull\/}
$\mix(A)$ of $ A$ consists of
$x\in {\mathbb V}^{(\mathbb B)}$
admitting an internal
family
$(a_{\xi})_{\xi \in \Xi}$
of elements
of $A$ and an internal partition
$(b _{\xi})_{\xi \in \Xi}$
of unity in $B$ such that $x$ is the mixing of
$(a_{\xi})_{\xi \in \Xi}$
by
$(b_{\xi})_{\xi\in\Xi}$;
i.e.,
$b_{\xi}x=b_{\xi}a_{\xi}$
for $\xi\in\Xi$ or, equivalently,
$x=\mix_{\xi\in\Xi}(b_\xi a_{\xi})$.

 Given a filter ${\mathscr F}$ on ${X\!\!\downarrow}$, let
$$
{{\mathscr F}\!\!\uparrow\downarrow}\!:=
\fil{({{\{F\!\!\uparrow\downarrow} \mid  F \in {\mathscr F}\}})}.
$$
Then
$\mix ({\mu ({\mathscr F})})=\mu ({{\mathscr F}\!\!\uparrow\downarrow})$
and
${{\mathscr F}\!\!\uparrow\downarrow}$
is the greatest cyclic filter coarser than ${\mathscr F}$.

The monad of ${\mathscr F}$ is called
{\it cyclic\/}
if $\mu ({\mathscr F})=\mix({\mu ({\mathscr F})})$.
Unfortunately, the cyclicity of a~monad
is not completely responsible for extensionality of a~filter.

%\section*{Monad Hulls}
The {\it  cyclic monad hull\/}
${\mu}_c (U)$
of an external set $U$ is defined as follows:
$$
x\in{\mu}_c (U)\,\leftrightarrow\, ({\forall}^{\st} V
=V\!\!\uparrow\downarrow ) V \supset U
\,\rightarrow\,x\in\mu(V).
$$
 If
$\mathbb B=\mathbb 2$,
then $\mu _c (U)$ is the monad of the standardization
of the external filter of supersets of $U$, i.e. the
({\it discrete$)$ monad hull\/} ${\mu}_d (U)$.

{\sl The cyclic monad hull of a~set is the cyclic
hull of its monad hull}
$$
{\mu}_c (U)=\mix ({{\mu}_d (U)}).
$$

%\section*{Essential Points}
A special role is played by
the {\it essential points\/}
of ${X\!\!\downarrow}$
constituting the external set
${}^e X$.
By definition, an essential point of
${}^e X $ belongs to the monad of some proultrafilter
on ${X\!\!\downarrow}$. The collection ${}^e X$ of all
essential points of $X$ is
usually external.

{\bf Test for Essentiality.}
{\it A point $x\in {}^e X$ if
and only if $x$ can be separated by a~standard set from every standard
cyclic set not containing $x$.}

If there is an~essential point
in the monad of an ultrafilter ${\mathscr F}$  then
$\mu  ({\mathscr F}) \subset  {}^e X $;
moreover,
${{\mathscr F} \!\!\uparrow\downarrow}$
is a~proultrafilter.

{\sl
A~filter $\mathscr F$ is extensional if and only
if   $\mu(\mathscr F) = \mu_c({}^e \mu{F})$.
A standard set $A$ is cyclic if and only if $A$ is the cyclic
monad hull of ${}^e A$}.

{\bf Test for the Mixing of Filters.}
{\it   Let
$({\mathscr F}_{\xi})_{\xi \in \Xi}$
be a~standard
family of extensional filters, and let
$(b_{\xi})_{\xi\in \Xi}$
be a~standard partition of unity. The filter
${\mathscr F}$ is the mixing
of
$({\mathscr F} _{\xi}) _{\xi \in \Xi}$
by $(b _{\xi }) _{\xi \in \Xi}$ if and only if
$$
(\forall ^{\St} \xi \in \Xi)
b _{\xi}\mu  ({\mathscr F}) =b _{\xi}\mu  ({\mathscr F} _{\xi}) .
$$
}

{\bf Properties of Essential Points.}
{\it
(1) The image of an essential point under an
extensional mapping is an essential point of the image;

(2) Let $E$ be a~standard set, and let $X$ be a~standard
element of
$ {\mathbb V}^{(\mathbb B)} $.
Consider the product
$X^{E^{\scriptscriptstyle\wedge}}$
inside
$ {\mathbb V}^{(\mathbb B)}$,
where
$E ^{\scriptscriptstyle\wedge}$
is the standard name of $E$
in
$\mathbb V^{(\mathbb B)} $.
If $x$ is an essential point of
${X^{E^{\scriptscriptstyle\wedge}}\!\!\downarrow}$
then for every
standard $e \in E $ the point
${x\!\!\downarrow}(e) $  is essential in
${X\!\!\downarrow}$;

(3)
Let ${\mathscr F}$ be a~cyclic filter in ${X\!\!\downarrow}$,
and let
${}^e\mu({\mathscr F})\!:= \mu(\mathscr F)\cap{}^e X$
be the set of essential points of
its monad. Then
$
{}^e\mu({\mathscr F}) ={}^ e \mu({\mathscr F}^{\uparrow\downarrow}).
$
}

%\section*{Procompactness}
Let $(X, {\mathscr U}) $ be a~uniform space inside
$ {\mathbb V}^{(\mathbb B)} $.The descent
$(X\!\!\downarrow , {\mathscr U}^\downarrow)$
is  {\it procompact\/}
or {\it cyclically compact\/}
if $(X, {\mathscr U})$ is compact inside
$ {\mathbb V}^{(\mathbb B)}$.
A similar sense resides in the notion of
{\it pro-total-boundedness\/}
and so on.

Every essential point of ${X\!\!\downarrow}$ is nearstandard,
i.e., infinitesimally close to a~standard point, if and only if
${X\!\!\downarrow}$
is procompact.

Existence of many procompact but not compact spaces provides a lot of
examples of inessential points.

{\bf Test for Proprecompactness.}
 {\sl A~standard space is the descent of a~totally bounded uniform space if
and only if its every essential point is prenearstandard, i.e. belongs
to the monad of a~Cauchy filter.}

 Let $Y$ to be a
universally complete vector lattice.
By Gordon's Theorem,  $Y={\mathscr R\!\!\downarrow}$ of the reals $\mathscr R$  inside
${\mathbb V}^{(\mathbb B)}$
over the base $\mathbb B:={\mathbb B}(Y)$ of~$Y$.

Denote by
$\mathscr E$
the filter of order units in $Y$, i.e.
${\mathscr E}\!:= \{ \varepsilon \in Y _+ \mid [\![\,\varepsilon =0\,]\!] ={\mathbb 0} \}$.

Put
$x \approx y \,\leftrightarrow\,(\forall^{\st}\varepsilon
\in {\mathscr E})\, (|x - y | < \varepsilon )$.
Given $a,b\in Y$,  write $a < b $
if
$[\![\,a < b\,]\!] ={\mathbb 1} $;
in other words,
$a>b\,\leftrightarrow\, a-b\in {\mathscr E}$.
Thus, there is some deviation from the understanding of the theory of ordered vector
spaces. Clearly, this is done in order to
adhere to the principles of introducing notations while descending and
ascending.

Let
${}^\approx Y$
be the {\it nearstandard part\/}
of $Y$.
Given
$y \in {}^ \approx Y $,
denote by ${}^\circ y $
(or by $\st (y) $) the {\it standard part\/}
of $y$, i.e. the unique standard element infinitely close to~$y$.

{\bf Tests for Order Convergence.}
{\it
For a~standard filter ${\mathscr F}$
in $Y$ and a~standard $z\in Y$, the following  are true:

(1)
$\inf_{F \in\mathscr F}\sup F\le z \,\leftrightarrow\,
(\forall y \in {}^{\bf .}
\mu(\mathscr F\!\!\uparrow\downarrow))\,
{}^\circ y \le z \,%\hfill\break\phantom{qqqqq}
\leftrightarrow(\forall y \in {}^e\mu({\mathscr F}\!\!\uparrow\downarrow))\,
{}^\circ y  \le z $;

(2)
$\sup_{F \in\mathscr F}\inf F\ge z\,\leftrightarrow\,
(\forall y\in {}^{\bf .}\mu(\mathscr F\!\!\uparrow\downarrow))\,
{}^\circ y\ge z\,\leftrightarrow
(\forall y\in {}^e\mu({\mathscr F}\!\!\uparrow\downarrow))\,
{}^\circ y\ge z$;

(3)
$\inf_{F \in\mathscr F} \sup F\ge z\,\leftrightarrow\,
(\exists y\in{}^{\bf .}\mu(\mathscr F\!\!\uparrow\downarrow))\,
{}^\circ y\ge z\,
%\hfill\break
%\phantom{qqqqq}
\leftrightarrow
(\exists y\in {}^e \mu({\mathscr F}\!\!\uparrow\downarrow))\,
{}^\circ y\ge z$;

(4)
$\sup_{F \in\mathscr F}\inf F\le z\,\leftrightarrow\,
(\exists y\in {}^{\bf .}\mu(\mathscr F\!\!\uparrow\downarrow))\,
{}^\circ y\le z\,
\leftrightarrow
(\exists y\in {}^e\mu(\mathscr F\!\!\uparrow\downarrow))
{}^\circ y \le z$;

(5)
${\mathscr F} \overset {(o)}\to
 z \leftrightarrow  ({ \forall y \in {}
^ e \mu  ({\mathscr F}\!\!\uparrow  \downarrow )})
 y  \approx z
 \leftrightarrow({\forall y\in\mu
({\mathscr F} ^{\uparrow \downarrow})}) y
 \approx  z $.
\smallskip

\noindent
Here
${}^{\bf .}\mu({\mathscr F}\!\!\uparrow
\downarrow)\!:= \mu ({\mathscr F}\!\!\uparrow \downarrow)
\cap {}^ \approx Y$,
and, as usual,
${}^e \mu ({\mathscr F}\!\!\uparrow \downarrow )$
is the set
of essential points of the monad
$\mu({\mathscr F}\!\!\uparrow\downarrow )$,
i.e.
${}^e\mu({\mathscr F}\!\!\uparrow\downarrow)
=\mu ({{\mathscr F}\!\!\uparrow\downarrow})\cap {}^e\mathscr R$.
}

\section*{Boolean Valued Monads}
Let us follow  the classical
approach of Robinson inside
${\mathbb V}^{(\mathbb B)}$.
In other words,  the classical and internal
universes and the corresponding $*$-map (Robinson's standardization)
are understood to be members of
${\mathbb V}^{(\mathbb B)}$.
Moreover, the nonstandard world is supposed  to be properly saturated.

The descent of the $*$-map is referred to as
{\it  descent standardization}.
Alongside the term  ``descent standardization''
the expressions like ``$B$-standar\-dization,'' ``prostandardization,''
etc. are in common parlance. Furthermore,
Denote the Robinson standardization of a~$B$-set $A$
by ${}^*\!A$.

The {\it descent standardization\/} of a~set
$A$ with $B$-structure, i.e. a~subset
of
${\mathbb V}^{(\mathbb B)}$,
is defined as
$({}^* (A\!\!\uparrow))\!\!\downarrow$
and is denoted
by ${}_* A$
(it is meant here that
${A\!\!\uparrow}$
is an element of the standard universe located inside
${\mathbb V}^{(\mathbb B)}$).

Thus,
${}^* a\in{}_* A\, \leftrightarrow\, a\in A\!\!\uparrow\downarrow$.
The {\it descent standardization\/}
${}_ * \Phi$
of an {\it extensional correspondence\/} $\Phi$
is also defined in a~natural
way.

Considering the descent standardizations of the standard names
of elements of the von Neumann universe ${\mathbb V}$,
use the abbreviations
${}^ *x\!:= {}^* (x ^{\scriptscriptstyle\wedge})$
and
${}_* x\!:= ({}^* x)\!\!\downarrow$
for
$x \in {\mathbb V} $.
The rules of placing and omitting  asterisks (by default) in
descent standardization are also assumed  as liberal as
those for the Robinson $*$-map.

{\bf Transfer.}
{\it
Let $\varphi =\varphi(x, y) $ be a~formula of  $\ZFC$
without any free variables other than $x$ and $y$.
Then
$$
(\exists x\in {}_* F)\,
[\![\,\varphi (x, {}^* z)\,]\!]={\mathbb 1}\,\leftrightarrow\,
(\exists x \in F\!\!\downarrow )\,
[\![\,\varphi (x, z)\,]\!] ={\mathbb 1};
$$
$$
(\forall x \in {}_ * F)\,[\![\,\varphi (x,{}^* z)\,]\!]={\mathbb 1}\,
\leftrightarrow\,(\forall x \in F\!\!\downarrow)\,
[\![\,\varphi (x, z)\,]\!]={\mathbb 1}
$$
for a~nonempty
element $F$ in
$ {\mathbb V}^{(\mathbb B)}$
and for every $z$.
}

{\bf Idealization.}
{\it
Let $X \!\!\uparrow$
and $Y$ be classical elements of
${\mathbb V}^{(\mathbb B)}$,
and
let $\varphi  =\varphi (x, y, z) $ be a~formula of
$\ZFC$. Then
$$
(\forall^{\fin}A\subset X)\,(\exists y\in {}_* Y)\,(\forall x\in A)\,
[\![\,\varphi ({}^* x, y, z)\,]\!]={\mathbb 1}
%\leftrightarrow\,
$$
$$
\leftrightarrow
(\exists y \in {}_* Y)\,
(\forall x \in X)\,[\![\,\varphi ({}^* x, y, z)\,]\!] ={\mathbb 1}
$$
for an internal element $z$ in
${\mathbb V}^{(\mathbb B)}$.
}

Given a~filter ${\mathscr F}$ of sets with $B$-structure, define the
{\it descent monad\/}
$m ({\mathscr F}) $ of~$\mathscr F$ as
$$
m ({\mathscr F})\!:=\bigcap\limits_{F \in\mathscr F}{}_* F.
$$

{\bf Meets of Descent Monads.}
{\it Let ${\mathscr E}$ be a~set of filters,
and let
${{\mathscr E}^\uparrow\!}:=\{\mathscr F^\uparrow \mid\mathscr F\in\mathscr E\}$
be its ascent to
$\mathbb V^{(\mathbb B)}$.
The following  are equivalent:

(1)
the set of  cyclic hulls  $\mathscr E$, i.e.
${{\mathscr E}\!\!\uparrow\downarrow\!}:=
\{{\mathscr F\!\!\uparrow\downarrow}\mid\mathscr F\in\mathscr E\}$,
is bounded above;

(2)
${\mathscr E}^\uparrow$
is bounded above inside
$\mathbb V^{(\mathbb B)}$;

(3)
$\bigcap\{ m(\mathscr F)\mid\mathscr F\in\mathscr E\}\ne\varnothing$.

 Moreover, in this event
$$
m(\sup\mathscr E\!\!\uparrow\downarrow)
=\bigcap\{ m(\mathscr F) \mid\mathscr F\in\mathscr E\};\quad
\sup {\mathscr E}^{\uparrow}=(\sup{\mathscr E})^{\uparrow}.
$$
}

It is worth noting that for an infinite set of descent monads, its
union, and even the cyclic hull of this union, is not a~descent monad
in general. The situation here is the same as for ordinary monads.

{\bf Nonstandard Tests for a~Proultrafilter}
{\it The following  are equivalent:

(1)
${\mathscr U}$ is a~proultrafilter;

(2)
${\mathscr U}$
 is an extensional filter with inclusion-minimal descent monad;

(3)
the representation
$\mathscr U =(x)^\downarrow\!:=
\fil{(\{{U\!\!\uparrow\downarrow} \mid x\in {}_* A \})}$
holds for each point $x$ of the descent monad
$m({\mathscr U});$

(4)
${\mathscr U}$
is an extensional filter whose
descent monad is easily caught by a~cyclic set; i.e.
either
$m(\mathscr U)\subset {}_* U$
or
$m(\mathscr U)\subset {}_*(X \setminus U)$
for every
$U =U\!\!\uparrow\downarrow$;

(5)
${\mathscr U}$
is a~cyclic filter satisfying the condition:
for every cyclic
$U$,
if
${}_* U\cap m(\mathscr A)\ne\varnothing$
then
$U\in\mathscr U$.
}

{\bf Nonstandard Test for the Mixing of Filters.}
{\it Let $({\mathscr F} _{\xi }) _{\xi \in \Xi }$ be
a~family of filters, let $(b _{\xi }) _{\xi \in \Xi}$
be a~partition of unity, and let
$\mathscr F =\mix_{\xi\in \Xi}(b_\xi\mathscr F_{\xi}^\uparrow)$
be the mixing of
$\mathscr F_\xi^\uparrow$
by
$b_{\xi}$.
Then
$$
m(\mathscr F^\downarrow)=\mix_{\xi\in \Xi}(b_\xi m(\mathscr F_\xi)) .
$$
}

A point $y$ of  ${}_ * X$ is called
{\it descent-nearstandard\/}
or simply {\it nearstandard\/}
if there is no danger of misunderstanding
whenever
${}^* x\approx y$
for some $x \in {X\!\!\downarrow}$;
i.e.,
$(x, y)\in m({\mathscr U}^\downarrow)$,
with ${\mathscr U}$  the uniformity on~$X$.

{\bf Nonstandard Test for Procompactness.}
{\it A set ${A\!\!\uparrow\downarrow}$ is procompact if and
only if every point of ${}_* A $ is descent-nearstandard.}

{\bf Truth Value on a~Proultrafilter.}
{\it Let $\varphi =\varphi (x) $ be a~formula of
$\ZFC$. The truth value of $\varphi $ is constant on the
descent monad of every proultrafilter ${\mathscr A}$; i.e.,
$$
(\forall x, y\in m (\mathscr A))\,
[\![\, \varphi (x)\,]\!]=[\![\,\varphi (y)\,]\!].
$$
}
Let
$\varphi =\varphi (x, y, z)$
be a~formula of $\ZFC$, and let
${\mathscr F}$ and   ${\mathscr G}$ be filters of sets with $B$-structure.

{\bf Rules of Descent Standardization.}
{\it The following
quantification rules are valid $($for internal $y$,  $z$ in
${\mathbb V}^{(\mathbb B)})$:

(1)
$(\exists x\in m(\mathscr F))\,
[\![\,\varphi (x, y, z)\,]\!]={\mathbb 1}
\leftrightarrow\,
(\forall F\in\mathscr F)\,(\exists x\in{}^*F)\,
[\![\,\varphi (x,y,z)\,]\!]={\mathbb 1}$;

(2)
$(\forall x\in m (\mathscr F))\,
[\![\,\varphi (x, y, z)\,]\!]={\mathbb 1}
\leftrightarrow\,(\exists F\in{\mathscr F}^{\uparrow\downarrow})
(\forall x \in {}_* F)\,[\![\,\varphi (x,y,z)\,]\!] ={\mathbb 1}$;

(3)
$(\forall x \in m (\mathscr F))\,(\exists y \in m (\mathscr G))
[\![\,\varphi (x,y,z)\,]\!]={\mathbb 1}$\hfill\break
 $\phantom{qqqqq}\leftrightarrow \,(\forall G\in\mathscr G)\,
(\exists F\in\mathscr F^{\uparrow\downarrow})\,
(\forall x\in {}^* F)\,(\exists y\in {}^* G)\,
[\![\,\varphi (x,y,z)\,]\!] ={\mathbb 1}$;

(4)
$(\exists x\in m (\mathscr F))\,(\forall y\in m (\mathscr G))\,
[\![\,\varphi (x,y,z)\,]\!] ={\mathbb 1}$\hfill\break
$\phantom{qqqqq}\leftrightarrow\,(\exists G \in\mathscr G^{\uparrow\downarrow})\,
(\forall F\in\mathscr F)\,(\exists x \in {}^* F)\,(\forall y\in {}^* G)\,
[\![\,\varphi (x, y, z)\,]\!]={\mathbb 1}$.

(5)
$({ \exists x \in m ({\mathscr F})}) [\![\, \varphi (x, {}^ * y, {}^ * z)\,]\!] ={\mathbb 1}
\leftrightarrow  (\forall F \in {\mathscr
F}) (\exists x \in F\!\!\uparrow  \downarrow )
  [\![\, \varphi (x, y, z)\,]\!] ={\mathbb 1};$

(6)
$({ \forall x \in m (
{\mathscr F})}) [\![\, \varphi (x, {}^ * y, {}
^ * z)\,]\!] ={\mathbb 1}$
$\leftrightarrow  (\exists F \in {\mathscr
F}^{\uparrow  \downarrow}) (\forall x \in F
) [\![\, \varphi (x, y, z)\,]\!] ={\mathbb 1};$

(7)
$({ \forall x \in m (
{\mathscr F})}) ({ \exists y \in m ({\mathscr
G})}) [\![\, \varphi (x, y, {}^ * z)\,]\!] ={\mathbb 1}$\hfill\break
$\phantom{qqqqq}\leftrightarrow  (\forall G \in {\mathscr
G}) (\exists F \in {\mathscr F} ^{\uparrow
  \downarrow}) (\forall x \in F) (\exists y
\in G\!\!\uparrow  \downarrow ) [\![\, \varphi (x,  y,  z)\,]\!] ={\mathbb 1}$;

(8)
$({ \exists x \in m (
{\mathscr F})}) ({ \forall y \in m ({\mathscr
G})}) [\![\, \varphi (x, y, {}^ * z)\,]\!] ={\mathbb 1}$\hfill\break
$\phantom{qqqqq} \leftrightarrow  (\exists G \in {\mathscr
G} ^{\uparrow  \downarrow}) (\forall F \in {\mathscr
F}) (\exists x \in F\!\!\uparrow  \downarrow )
(\forall y \in G)  [\![\, \varphi (x,  y,  z)\,]\!] ={\mathbb 1}.$
}

\section*{The Escher Rules in Vector Lattices}

The fact that $E$ is a~vector lattice is
a~restricted formula, say, $\varphi (E, {\mathbb R}) $. Hence,
recalling the restricted transfer principle, we come to the
equality
$[\![\,\varphi (E^{\scriptscriptstyle\wedge},
{\mathbb R}^{\scriptscriptstyle\wedge})\,]\!] ={\mathbb 1}$;
i.e.,
$E^{\scriptscriptstyle\wedge}$
is a~vector
lattice over the ordered field
${\mathbb R}^{\scriptscriptstyle\wedge}$
inside
$ {\mathbb V}^{(\mathbb B)}$.

Let
$E^{\scriptscriptstyle\wedge\sim}$~be the space of
regular
${\mathbb R}^{\scriptscriptstyle\wedge}$-linear
functionals from
$E^{\scriptscriptstyle\wedge}$
to $\mathscr R$.
It is easy  that
$E^{\scriptscriptstyle\wedge\sim}\!:=
L^\sim (E^{\scriptscriptstyle\wedge}, \mathscr R)$
is a~$K$-space, i.e. a Dedekind complete vector lattice,   inside
${\mathbb V}^{(\mathbb B)}$.
Since
$E^{\scriptscriptstyle\wedge\sim}$
is a~$K$-space,
the descent
${E^{\scriptscriptstyle\wedge\sim}\!\!\downarrow}$
of
$E^{\scriptscriptstyle\wedge\sim}$
is a~$K$-space too.

Turn to the universally complete vector lattice
$F\!:= {\mathscr R\!\!\downarrow}$. For every operator
$T \in L^ \sim (E, F) $
the ascent ${T\!\!\uparrow}$ is defined by
the equality
$[\![\, Tx ={T\!\!\uparrow}(x^{\scriptscriptstyle\wedge})\,]\!] ={\mathbb 1}$
for all $x \in E$.
If
$\tau\in E^{\scriptscriptstyle\wedge\sim}$,
then
$[\![\, \tau :  E^{\scriptscriptstyle\wedge} \to \mathscr R \,]\!] ={\mathbb 1} $;
hence,
the operator
${\tau \!\!\downarrow}: E \to F $
is available.
Moreover,
${\tau\!\!\downarrow\uparrow}=\tau $.
On the
other hand,
${T\!\!\uparrow\downarrow} =T $.

{\sl  For every
$T\in L^{\sim}(E, F)$
the ascent
${T\!\!\uparrow}$
is a~regular
${\mathbb R}^{\scriptscriptstyle\wedge}$-functional on
$E^{\scriptscriptstyle\wedge}$
inside
${\mathbb V}^{(\mathbb B)}$;
i.e.,
$[\![\,{T\!\!\uparrow}\in E^{\scriptscriptstyle\wedge\sim}\,]\!] ={\mathbb 1}$.
The mapping
$T\mapsto {T\!\!\uparrow}$
is a~linear and lattice isomorphism  between
$L^{\sim}(E, F)$ and
${E^{\scriptscriptstyle\wedge\sim}\!\!\downarrow}$.}

An operator
$S\in L^\sim (E, F) $
is  a~{\it fragment\/} or {\it component\/}
of $0\leq T\in L^\sim (E, F)$ if
$S \wedge (T - S)=0$.
Say that  $T$ is
{\it $F$-discrete\/}
whenever
${[0, T ]=[0, I_F] \circ T}$;
i.e.,
for every
$0 \le S \le T$
there is an operator
$0 \le \alpha \le I_F$
satisfying
${S=\alpha\circ T}$.
Let
$L_a^\sim (E, F)$
be the band of
$L^\sim (E, F)$
generated by $F$-discrete
operators, and write
$L_d^\sim (E, F)\!:= L _a^\sim (E, F)^\perp $.
The bands
$(E^{\scriptscriptstyle\wedge\sim})_a $
and
$(E^{\scriptscriptstyle\wedge\sim}) _d$
are introduced similarly. The elements of
$L_d^\sim (E, F) $
are usually referred to as {\it $F$-diffuse\/}
operators.
The ${\mathbb R}$-discrete or ${\mathbb R}$-diffuse operators
are called for the sake of brevity {\it discrete\/} or {\it diffuse\/}
functionals.

{\bf Rules of Descending.}
{\it
Consider
$S, T \in L^\sim (E, F) $
and put  $\tau\!:= T\!\!\uparrow$;
$\sigma\!:= S\!\!\uparrow$.
The following  are true:

(1)
$T\ge 0\, \leftrightarrow\,
[\![\,\tau \ge 0\,]\!] ={\mathbb 1}$;

(2)
$[\![\,S$ is a~fragment of $T\,]\!]\,\leftrightarrow\,[\![\,\sigma$
is a~fragment of
$\tau\,]\!]={\mathbb 1}$;

(3)
$[\![\,T$ is $F$-discrete\,$]\!]\leftrightarrow[\![\,\tau
\text{ is discrete}\,]\!]= {\mathbb 1}$;

(4)
$T \in L_a^\sim (E, F)\,\leftrightarrow\,[\![\,\tau
\in (E ^{\scriptscriptstyle\wedge\sim})_a\,]\!]
={\mathbb 1} $;

(5)
$T \in L_d^\sim (E, F)\,\leftrightarrow\,[\![\,\tau
\in (E^{\scriptscriptstyle\wedge\sim})_d\,]\!]={\mathbb 1}$.

(6)
$[\![\,T$ is a lattice homomorphism\,$]\!]$
$\,\leftrightarrow\,[\![\,\tau$ is a lattice homomorphism$\,]\!] ={\mathbb 1} $.

}

Let $E$ stand for a~vector lattice and $F$,
for a~$K$-space.
A set ${\mathscr P}$ of band projections in
$L^{\sim}(E, F)$
{\it generates the fragments of $T$},
$0\leq T\in L^\sim (E, F)$, provided that
$Tx^+ =\sup \{pTx \mid p \in {\mathscr P}\}$
for all $x \in E$.
If this happens for
all
$0\le T\in L^{\sim}(E, F)$,
 then
${\mathscr P}$ is a~{\it generating} set.

Put
$F\!:= {\mathscr R\!\!\downarrow}$
and let $p$ be a~band projection in $L^\sim (E, F)$. Then
there is a~unique element
${p\!\!\uparrow}\in {\mathbb V}^{(\mathbb B)}$
such that
$[\![\,\, p\!\!\uparrow$ is a~band projection in
$E^{\scriptscriptstyle\wedge\sim}\,]\!]={\mathbb 1}$
and
$(p T)\!\!\uparrow =p\!\!\uparrow T\!\!\uparrow $
for all
$T\in L^\sim (E, F)$.

{\bf Rules of Fragmenting.}
{\it
Consider some set ${\mathscr P}$ of band projections in
$L^\sim (E, F)$ and a~positive operator
$T\in L^\sim (E, F)$. Put
$\tau := {T\!\!\uparrow}$ and ${{\mathscr P}\!\!\uparrow}:=
\{{p\!\!\uparrow}\mid p \in {{\mathscr P} \}\!\!\uparrow}$.
Then $[\![\,{{\mathscr P}\!\!\uparrow}$
is a~set of  band projections
in~$E^{\scriptscriptstyle\wedge\sim}\,]\!] ={\mathbb 1} $
and the following
are true:

(1)
$[\![\,{\mathscr P}$ generates the fragments of $T\,]\!]\leftrightarrow\hfill\break
\phantom{qqqqqq}[\![\,{\mathscr P}\!\!\uparrow$
generates the fragments of $\tau\,]\!]=\mathbb 1;$

(2)
$[\![\,{\mathscr P}$  is a~generating set\,$]\!]$
$\leftrightarrow\,[\![\,{\mathscr P}\!\!\uparrow$
 is a~generating set$\,]\!]={\mathbb 1}$.
}

Given a~set $A$ in a~$K$-space,  denote by
$A^{\vee}$
the result of adjoining to $A$ suprema of  every
nonempty finite subset of~$A$. Let
$A^{\uparrow}$
stand for the result of adjoining to $A$ suprema of nonempty increasing
nets of elements of~$A$.
The symbols
$A^{\uparrow\downarrow}$
and
$A^{\uparrow\downarrow\uparrow}$
are understood naturally (cp.~\cite{Pagter}--\cite{KusKut2005}).

Put ${\mathscr P}(f)\!:= \{ pf \mid p\in {\mathscr P}\}$
and note that  $E$ will for a time being stand for
a~vector lattice over a~dense subfield of ${\mathbb R}$ while
${\mathscr P}$ is a~set of band projections in $E^{\sim}$.
Let ${\mathfrak E}(f)$ stand  for the set of all fragments of~$f$.

{\bf Up-Down Theorem.}
{\it
The following  are equivalent:

(1)
${\mathscr P(f)}^{\vee(\uparrow \downarrow \uparrow)}={\mathfrak E}(f)$;

(2)
${\mathscr P}$ generates the fragments of $f$;

(3)
$(\forall x \in {}^{\circ}E)
(\exists p \in {\mathscr P}) pf(x)\approx f(x^{+}) $;

(4)
a~functional $g$ in
$[0, f]$
is a~fragment of
$f$
if and only if
$$
\inf\limits_{p \in {\mathscr P}}({p^\perp g(x) +p(f-g)(x)})=0
$$
for every
$0\leq x\in E$;

(5)
$(\forall g \in {}^{\circ}{\mathfrak E}(f))
(\forall x \in {}^{\circ}E_{+})(\exists p \in {\mathscr P})
| p f -g |(x)\approx 0$;

(6)
$\inf\{|pf-g| (x) \mid p \in {\mathscr P}\}=0$
for all fragments
$g \in {\mathfrak E} (f)$
and $x\ge 0$;

(7)
for $x \in E_{+}$
and
$g \in {\mathscr E}(f)$
there is an element
$p\in{\mathscr P}(f)^{\vee(\uparrow\downarrow\uparrow)}$,
satisfying
$$|pf-g|(x) =0.$$}

{\scshape Proof.}
The implications
(1)~$\rightarrow$~(2)~$\rightarrow$~(3) are obvious.

(3)~$\rightarrow$~(4):
We will work within the {\it standard entourage};
i.e., we presume that all free variables are standard.
Note first  that validity of the sought equality for all
functionals $g$ and $f$ satisfying $0\leq g\leq f $
amounts to
existence of
$p\in {\mathscr P}$,
given a~standard $x \ge  0$,
such that
$p^\perp g (x) \approx 0$
and
$p(f-g)(x)\approx 0$.
(As usual, $p^\perp$ is the {\it complementary band projection\/}
to $p$.)
Thus,
${}^ \circ p ({g \wedge (f-g)}) (x) \leq  {}^\circ p (f-g) (x) =0$
and
${}^{\circ}p^{\perp}({(f-g)\wedge g})(x)\leq {}^ \circ p^\perp g(x)=0$,
i.e.
$g\wedge (f-g)=0$.

Prove now that, on assuming (3), the sought equality
ensues from the conventional criterion for disjointness:
$$
\inf \{ g (x_1) + (f-g) (x_2) \mid x_1 \ge 0,\, x_2 \ge 0,\, x_1 + x_2 =x \}=0.
$$

Given a~standard $x$, find internal positive $x_ 1 $ and
$x_ 2 $ such that $x =x_ 1  + x_ 2 $ and, moreover,
$g (x_1) \approx 0$
and
$f (x_2) \approx g(x_2)$.
By  (3), it follows from the Kre\u\i n--Milman Theorem
that the fragment $g$ belongs to the weak
closure of ${\mathscr P} (f) $. In particular, there is an element
$p \in {\mathscr P} $
satisfying
$g (x_1) \approx pf(x_1)$
and
$g (x_2) \approx pf(x_2)$.
Thus,
$p^\perp g(x_2) \approx 0$,
because
$p^\perp g\leq p^\perp f$.
Finally,
$p^\perp g (x) \approx 0$.
Hence,
$$
\gathered
p(f-g)(x) =pf(x_2) + pf(x_1)-pg (x)
\approx\\
g(x_2)+
g(x_1)-pg(x)\approx p^\perp g (x) \approx 0.
\endgathered
$$
This yields the claim.

(4)~$\rightarrow$~(5):
Using the equality
$|pf-g|(x) =p^\perp g(x)+ p(f-g) (x)$,
we may find
$p\in{\mathscr P} $
so that
$p^\perp g(x)\approx 0$
and
$p(f-g)(x)\approx 0$.
This justifies the claim.

The equivalence (5)~$\leftrightarrow$~(6) is clear.
The implications (5)~$\rightarrow$~(7)~$\rightarrow$~(1)
are standard.
The proof is complete.

We now turn to  principal bands.
For positive functionals $f$
and $g$ and for a~generating set of band projections
${\mathscr P}$, the following  are equivalent:

(1)
$g \in \{f\}^{\perp \perp}$;

(2)
If $x$ is a {\it limited} element of $E$, i.e.
$x \in {}^{\fin} E :=
\{ x\in E \mid (\exists \overline x \in {}^\circ E)
|x| \leq \overline x\}$,
then
$pg (x) \approx 0$
whenever
$pf (x) \approx 0$
for
$p \in {\mathscr P}$;

(3)
$(\forall x \in E_ +)
(\forall\varepsilon >0) (\exists \delta  >0)
(\forall p\in {\mathscr P}) pf (x) \leq  \delta  \to  pg (x) \leq \varepsilon$.

\noindent With the principal bands available, we may proceed  to the  principal projections.

 Let $f$ and  $g$ be positive
functionals on $E$, and let $x$ be a~positive element of $E$.
The following representations of the band projection  onto
Denote the band projection to
$\{f\}^{\perp\perp}$ by~$b_ f $

{\bf Principal Projection on a Functional.}
{\it The following representations hold:

(1)
$ b_f g(x)\Rightharpoonup\inf{}^\ast\{{}^\circ pg(x)\mid
p^\perp f(x)\approx 0,\,p\in \mathscr P\}$,\hfill\break
\phantom{ggggg} where $\Rightharpoonup$
means that the formula is {\it exact},
i.e.,  equality is attained;

(2)
$ b_fg(x)=\sup_{\varepsilon>0}\inf\{pg(x) \mid
p^\perp f(x)\leq\varepsilon,\,p \in {\mathscr P}\};$

(3)
$ b_ f g(x)\Rightharpoonup\inf {}^\ast\{{}^\circ g(y) \mid
f(x-y)\approx 0,\, 0\leq y \leq x\}$;

(4)
$(\forall\varepsilon>0)\,(\exists \delta >0)\,
(\forall p\in{\mathscr P})\,p f (x)<\delta\rightarrow b_fg (x)
\leq p^\perp g(x)+\varepsilon$;

\hskip20pt $(\forall\varepsilon>0)\,(\forall \delta>0)\,
(\exists p\in{\mathscr P})\,pf(x)<\delta\wedge p^\perp g(x)
\leq  b_fg (x)+\varepsilon$;

(5)$(\forall\varepsilon>0)\,(\exists\delta>0)\,
(\forall 0\leq y\leq x)\,f(x-y)\leq\delta\rightarrow
 b_f g(x)\leq g(y)+\varepsilon$;

\hskip20pt$(\forall\varepsilon>0)\,(\forall\delta>0)\,
(\exists 0\leq y\leq x)\,f(x-y)\leq\delta\wedge g(y)
\leq b_fg(x)+\varepsilon $.
}

Ascending to and descending from the appropriate Boolean valued universe,
we implement principal bands in the operator case.

{\it For a~set of band projections
${\mathscr P}$ in $L^{\sim}(E, F)$
and
$0\leq S\in L^{\sim} (E, F)$
the following  are equivalent:

(1)
$\mathscr P(S)^{\vee(\uparrow
\downarrow \uparrow)}={\mathfrak E} (S)$;

(2)
${\mathscr P}$ generates the fragments of $S$;

(3)
$T \in [0,S]$
is a~fragment of
$S$
if and only if
$$
\inf_{p\in\mathscr P}({p^\perp Tx + p(S -T)x})=0
$$
\phantom{ggggg}for all
$0\leq x \in E$;

(4)
$(\forall x\in {}^\circ E)\,
(\exists p\in{\mathscr P}\!\!\uparrow\downarrow)\,pSx\approx Sx^ +$.
}

Using the simplest Escher rules and Nelson's algorithm yields
the description of the~principal band generated by an operator.

{\it For positive operators $S$ and $T$ and a~generating
set ${\mathscr P}$ of band projections in $L^{\sim}(E, F)$, the
following  are equivalent:

(1)
$T\in\{S\}^{\perp\perp}$;

(2)
$(\forall x\in {}^{\fin} E)\,
(\forall p\in {\mathscr P})\,(\forall b\in \mathbb B)\,
 b pSx\approx 0\rightarrow  b pTx \approx 0$;

(3)
$(\forall x\in {}^{\fin} E)\,
(\forall b\in \mathbb B)\, b Sx\approx 0\rightarrow b Tx\approx 0$;

(4)
$(\forall x\ge 0)\,(\forall\varepsilon\in{\mathscr E})\,
(\exists\delta\in{\mathscr E})\,(\forall p\in {\mathscr P})\,(\forall b\in \mathbb B)\,
 b pSx\leq\delta\rightarrow b pTx\leq\varepsilon$;

(5)
$(\forall x\ge0)\,(\forall \varepsilon \in {\mathscr E})\,
(\exists \delta \in {\mathscr E})\,(\forall b\in \mathbb B)
 b Sx\leq\delta\rightarrow  b Tx\leq \varepsilon$.
}

 Let $E$ be a~vector lattice,
and let $F$ be a~$K$-space having  the filter of order units ${\mathscr E}$
and the base $\mathbb B$. Suppose that $S$ and $T$ are positive operators in
$L^{\sim}(E, F)$  and $R$ is the band projection of $T$ to
the band $\{S\}^{\perp \perp }$.

{\bf Theorem of Principal Projection.}
{\it For
a positive $x \in E $, the following  are valid:

(1)
$Rx=\sup_{\varepsilon\in\mathscr E}
\inf\{ b Ty +  b^\perp Sx \mid 0\leq y\leq x,
 b\in \mathbb B,\, b S(x-y)\leq\varepsilon\}$;

(2)
$Rx=\sup_{\varepsilon\in\mathscr E}
\inf\{( b p)^\perp Tx \mid  b pSx\leq\varepsilon,\,
p\in{\mathscr P}, b\in \mathbb B\}$,

\noindent
where ${\mathscr P}$ is a~generating set of band projections in~$F$.
}

In closing, turn to the revisited  Farkas Lemma (cp.\cite{Trends}, \cite{Farkas} and \cite{Polyhedral}).
Let $X$ be a $Y$-seminormed real vector space, with $Y$
a~$K$-space. Given are some dominated polyhedral sublinear
operators  $P_1,\dots, P_N$ from $X$ to $Y$ and a dominated
sublinear operator $P: X\to Y$.

{\bf Polyhedral Lagrange Principle.}
{\it The finite value of the
constrained problem
$$
P_1(x)\le u_1,\dots,P_N(x)\le u_N,\quad P(x)\to\inf
$$
is the value of the unconstrained problem for
an appropriate Lagrangian without any constraint
qualification but polyhedrality}.

Polyhedrality is omnipresent and so finds applications in inexact data processing (cp.~\cite{Fiedler}).
Let $X$ be a   $Y$-seminormed real space, with $Y$
a $K$-space. Assume given a dominated polyhedral
sublinear operator $P:X \to Y$,
a dominated sublinear operator $Q: X\to Y$, and
$u,v\in Y$. Assume further that  $\{P\le u\}\neq\varnothing$.

{\bf Interval Farkas Lemma.}
{\it The following are equivalent:

{\bf(1)} for all  $b\in \mathbb B$, with $\mathbb B$ the base of~$Y$,  the sublinear operator inequality
$bQ\circ\sim (x)\ge -bv$
is a~consequence of the polyhedral sublinear operator inequality
$bP(x)\le bu$,
i.e.,
$
\{bP\le bu \}\subset\{bQ\circ\sim\ge -bv\},
$
 with $\sim(x):=-x$ for all $x\in X$;

{\bf(2)} there are $A\in\partial(P)$, $B\in\partial(Q)$, and
a positive orthomorphism $\alpha\in\Orth(m(Y))$ on the universal completion
$m(Y)$ of~$Y$ satisfying
$
B=\alpha A,\quad
\alpha u\le v.
$
}

%\newpage

\bibliographystyle{plain}

%\newpage
%\theendnotes

\end{document}